\documentclass[12pt,leqno,fleqn,epsfig]{article}
\usepackage{amssymb, epsfig, amsmath, amsthm}
\usepackage{color}
\usepackage{mathrsfs}       

\newcommand{\R}{\mathbb{R}}

\newcommand{\Z}{\mathbb{Z}}

\newcommand{\N}{\mathbb{N}}

\newcommand{\supp}{\operatorname*{supp}}

\newcommand{\bb}{\begin{equation}}
\newcommand{\ee}{\end{equation}}
\newcommand{\bea}{\begin{eqnarray}}
\newcommand{\eea}{\end{eqnarray}}
\newcommand{\bean}{\begin{eqnarray*}}
\newcommand{\eean}{\end{eqnarray*}}

\newcommand{\var}{\varepsilon}

\newcommand{\intl}{\int\limits}

\newcommand{\Beweisende}{\rule{0.2cm}{0.2cm}}

\newcounter{secnum}

\newtheorem{thm}{Theorem}[section]

\newtheorem{lem}[thm]{Lemma}

\theoremstyle{definition}
\newcommand{\n}{\nonumber}
\newtheorem{defin}[thm]{Definition}
\newtheorem{rem}[thm]{Remark}

\renewcommand{\l}{\lambda}
\newcommand{\bq}{\begin{eqnarray}}
\newcommand{\eq}{\end{eqnarray}}
\newcommand{\bqn}{\begin{eqnarray*}}
\newcommand{\eqn}{\end{eqnarray*}}
\title{Removing  discretely self-similar singularities for the 3D Navier-Stokes equations} 
 
\author{Dongho Chae$^*$  and J\"{o}rg Wolf $^\dagger$\\
\ \\
 $*$Department of Mathematics\\
Chung-Ang University\\
 Seoul 156-756, Republic of Korea\\
 e-mail: dchae@cau.ac.kr\\
and \\
$\dagger$Department of Mathematics\\
Humboldt University Berlin\\
Unter den Linden 6, 10099 Berlin, Germany\\
e-mail: jwolf@math.hu-berlin.de}
\date{}
\begin{document}

\maketitle
\begin{abstract}
We study the scenario of discretely self-similar blow-up for Navier-Stokes equations.
We prove that  at the possible blow-up  time such solutions  only one point singularity. 
 In case of the scaling parameter  $ \lambda $ near $ 1$ we remove the singularity.
  \noindent{\bf
keywords:} Navier-Stokes equations, discretely self-similar singularities

\end{abstract}

\section{Introduction}
\label{sec:-1}
\setcounter{secnum}{\value{section} \setcounter{equation}{0}
\renewcommand{\theequation}{\mbox{\arabic{secnum}.\arabic{equation}}}}

We consider the Navier-Stokes equations in $ \R^{3}$
\begin{equation}
\partial _t u  + (u \cdot \nabla ) u -\Delta u  = - \nabla \pi,\qquad  \nabla \cdot  u =0,
\label{1.1}
\end{equation}
where $u= (u_1, u_2, u_3)$ denotes the velocity field, while $ \pi $ stands for the  pressure.  
We are concerned on the the (backward) self-similar type blow-up of the smooth solutions to the Navier-Stokes equations.
We say that a function $u: \Bbb R^3 \times (-\infty, 0) \to \Bbb R^3$ is self-similar of 
$ u(x,t)=\l u(\l x, \l^2 t)$ for all $\ >1$ and for all $(x,t) \in \Bbb R^3 \times (-\infty, 0) $.
A self-similar function has the representation
$u(x,t)=\frac{1}{\sqrt{-t}} U(\frac{x}{\sqrt{-t}})$ for a function $U: \Bbb R^3\to \Bbb R^3$, which is called the profile of $u$.
Leray first considered and  asked the question of possibility of the self-similar blow-up for the 3D Navier-Stokes equations in \cite{ler}.
Later, Ne\v{c}as, R\r{u}\v{c}i\v{z}ka and \v{S}ver\'ak proved in \cite{NeRuSv1996} that there exists no nontrivial solution to the  3D Navier-Stokes equations, having  the profile $U\in L^3$. This result was extended by Tsai in \cite{Tsai1998}, showing  the triviality of a self-similar solution to the Navier-Stokes equation,  
which satisfies the local energy inequality, or the profile $U$ of which  belongs to $ L^p(\R^{3})$ for some $ p\in (3, +\infty]$. 

\hspace{0.5cm}
For more general notion of the discretely self-similar solutions Tsai\cite{Tsai2014} proved existence of {\em forward} discretely self-similar solutions for the scaling parameter $\lambda$ close to one, while in more recent paper \cite{brad} Bradshaw and Tsai proved existence of the global {\em forward} discretely self-similar solutions to the Navier-Stokes equations for arbitrary $\lambda$.

Now it becomes quite  natural to ask whether such result also holds for  {\em backward} discretely self-similar solutions.  
Here we call $ u: \R^{3}\times (-\infty,0) \rightarrow \R^{3}$  backward discretely self-similar with respect to 
$ \lambda \in (1, +\infty)$ or shortly $ \lambda$-DSS if for all $ (x,t)\in \R^{3}\times (-\infty, 0)$ it holds 
\begin{equation}
u(x,t)=\lambda u(\lambda x, \lambda^2 t). 
\label{1.2}
\end{equation}
Defining $ u_\lambda (x,t) =\lambda u(\lambda x, \lambda^2 t)$, the relation \eqref{1.1} becomes $ u_\lambda =u$. We recall that the notion self-similarlity implies $u=u_\l$ for all $\l >1$,
while in the case of discrete self-similarity the scaling parameter $\l >1$ is a fixed number.
In the case of Euler equations the nonexistence results for the backward discretely self-similar solutions  are obtained in \cite{ChaeTsai2014, cha2}. For the case 
of Navier-Stokes equations such nonexistence result for the nontrivial-discretely self-similar solutions  is still not available in the literature 
(see Remark 1.2 below for the case $ u \in C((- \infty, 0); L^3(\R^{3}))$\,). As  stated in \cite[section 1]{Tsai2014} it is an open problem in the nontrivial profile case.

\hspace{0.5cm}
In what follows we set $ Q= \R^{3}\times (-\infty,0)$. The first aim of this paper is to prove that such backward $ \lambda $- DSS solutions  
to the Navier-Stokes equations in $ \overline{Q} $ are regular outside  the origin  $ z=(0,0)$, and they behave like $ \frac{C_{ \ast}}{\,\sqrt{-t}+| x|}$ for a positive constant $ C_{ \ast}>0$. Our second main theorem is the existence of $ \lambda _{ \ast} >1$ depending on $ C_{ \ast}$ such that for all $ \lambda \in (1, \lambda _{ \ast}]$ every  backward $ \lambda $-DSS solution $ u$ to the Navier-Stokes equations with $| u(x,t)| \le  \frac{C_{ \ast}}{\,\sqrt{-t}+| x|}$ is trivial. 

\hspace{0.5cm}
For $ z_0=(x_0, t_0)\in \R^{4}$ and $ 0<r<+\infty$ by $ Q(z_0, r)$ we denote  the parabolic cube $ B(x_0, r)\times (t_0-r^2, t_0)$. 
Here $ B(x_0, r)$ stands for the usual ball in $ \R^{3}$ with respect to the Euclidian norm. We set 
\[
\hspace*{-1cm}V^2(Q(z_0,R)) = L^\infty(-t_0-R^2, t_0; L^2(B(x_0,R))) \cap L^2(t_0-R^2, t_0; W^{1,\, 2}(B(x_0,R))). 
\]
For any Banach space $ X$ of vector functions  by $ X_\sigma $ we denote the subspace of all divergence free fields.  

 \vspace{0.3cm}
 \hspace{0.5cm}
 Our first main result shows that for each $\l >1$ the $ \lambda $-DSS solution is regular everywhere except at one point.
 \begin{thm}
 \label{thm1.3}
For  $ 3 \le  p < +\infty$  let $ u \in C((- \infty, 0); L^p(\R^{3}))\cap C^\infty(Q)$ be a  solution to the 
Navier-Stokes equations, and $ \lambda $-DSS for some $ \lambda \in (1, +\infty)$.  Then the solution $ u $ is regular on $ \overline{Q}   \setminus \{(0,0)\}$, and
satisfies the estimate
\begin{equation}
| u(x,t)| \le \frac{C}{\,\sqrt{-t}+| x|}\quad  \forall\,(x,t)\in Q. 
\label{1.6a}
\end{equation}
\end{thm}

\begin{rem}
If  $ u \in C((- \infty, 0); L^3(\R^{3}))$, and discretely self-similar,  then $ u\in L^\infty(-\infty, 0; L^3(\R^{3}))$.  Thus, in case $ p=3$, by using the result in \cite{esc}, we get the full regularity $ u$ in $ \overline{Q}  $. 
\end{rem}
 
\hspace{0.5cm}
Our second  main result of this paper is to show that for $\l>1$ close to $1$ one can remove the $ \lambda $-DSS solution.

\begin{thm}
\label{thm1.4}
For evrey $ C_{ \ast}>0$ there exists $ \lambda _{ \ast}>1$ depending on $ C_{ \ast}$ such that if 
 $ u \in C^\infty(Q)$ is  a $ \lambda$-DSS solution the Navier-Stokes equations for $ \lambda \in (1, \lambda _{ \ast})$, which satisfies 
\begin{equation}
| u(x,t)| \le \frac{C_{ \ast}}{\,\sqrt{-t}+| x|}\quad  \forall\,(x,t)\in Q. 
\label{1.7}
\end{equation}
Then $ u \equiv 0$. 
\end{thm}
 
 \begin{rem}
 Note that the criterion of \cite{gus} implies that if $ C_{ \ast}$  in \eqref{1.7}  is small enough, then every   $ \lambda $-DSS solution 
 to the Navier-Stokes equations satisfying  \eqref{1.7}  is trivial.  
 \end{rem}
 The notion of asymptotically self-similar scenario of solutions to semi-linear heat equations  has been introduced first by Giga and Kohn 
in \cite{giga}.    As an application of Theorem \ref{thm1.4} we can exclude a scenario of asymptotically discretely self-similar singularities with the scaling parameter $\l$ close to $1$.

  \begin{thm}
\label{thm1.5}
Let $3\leq p< +\infty$ and $ u \in C([0, t_*); L^p (\Bbb R^3))$ be a local in time smooth solution to (NS). Suppose there exists  $v(x,t)\in  C((-\infty, 0); L^p (\Bbb R^3)) $, fulfilling the inequality \eqref{1.7},  which is a $\l-$DSS 
function  with respect to $(x_*, t_*)$ with  $\l \in (1, \l_*)$ for $\l_*$ according to  Theorem \ref{thm1.4} such that 
\bb\label{as}
\lim_{t\to t_*} (t_*-t)^{\frac{p-3}{2p}} \sup_{t<\tau<t_*} \left\| u(\cdot, \tau)- v(\cdot, \tau)\right\|_{L^p (B_{R\sqrt{t_*-t}} (x_*)) } =0 \quad \forall R>0.
\ee
Then, $v=0$, and $(x_*, t_*)$ is a regular point.
\end{thm} 
 \begin{rem} We are assuming that $v$ in  (\ref{as}) is  a DSS function(not DSS solution of NS), and due to the factor $(t_*-t)^{\frac{p-3}{2p}} \to 0$  for $p>3$
 the ``convergence" $u\to v$ in $ L^p (B_{R\sqrt{t_*-t}} (x_*))$  is not guaranteed in general.
\end{rem}
\section{Regularity and decay  for $ \lambda$-DSS  solutions to the Navier-Stokes equations-Proof of Theorem\,\ref{thm1.3}}
\label{sec:-2}
\setcounter{secnum}{\value{section} \setcounter{equation}{0}
\renewcommand{\theequation}{\mbox{\arabic{secnum}.\arabic{equation}}}}

Let $ u\in C((- \infty, 0); L^p(\R^{3}))\cap C^\infty(Q)$ be a  solution to the Navier-Stokes equations, and 
$ \lambda $-DSS for some $ \lambda \in (1,+\infty)$. 

\vspace{0.2cm}
{\it \underline{1. Asymptotical behavior in time}}: We  prove that 
\begin{equation}
\| u (t)\|_p^p \le 
\lambda ^{ p-3}\| u \|^p_{ L^\infty(-\lambda ^2, -1; L^p) }(-t)^{ \frac{3-p}{2}}\quad \forall\,t\in (-\infty, 0). 
\label{2.4a}
\end{equation}

\hspace{0.5cm}
Let $ t\in (-\infty, 0)$ be arbitrarily chosen. Clearly, there exists a unique $ k\in \Z$ such that 
$ t\in  [-\lambda^{ 2(k+1)}, - \lambda ^{ 2k})$.  Recalling that $ u$ is $ \lambda$-DSS, we calculate    
\begin{align*}
\| u(\lambda^{ -2k}t ) \|_p^p&= \intl_{ \R^{3}} | u (x, \lambda ^{-2k }t)|^p  dx =    
\lambda ^{(p -3)k}\intl_{ \R^{3}} |\lambda ^{ -k} u (\lambda^{ -k} x, \lambda ^{-2k }t)|^p  dx
\\
&= \lambda ^{(p -3)k} \| u (t)\|_p^p.
\end{align*}
Since $ \lambda^{ -2k}t\in [ -\lambda ^2, -1)$ and $ \lambda ^{ 2k} < -t \le  \lambda ^{ 2(k+1)} $ we get 
\begin{align*}
\| u (t)\|_p^p &=\lambda ^{(3 -p)k}\| u(\lambda^{ -2k}t ) \|_p^p
\\
 &\le  \lambda ^{(3 -p)k} \| u \|^p_{ L^\infty(-\lambda ^2, -1; L^p) }
=
\lambda ^{ p-3} \| u \|^p_{ L^\infty(-\lambda ^2, -1; L^p) }(\lambda^{2 (k+1)}) ^{ \frac{3-p}{2}} 
\\
&\le  \lambda ^{ p-3}\| u \|^p_{ L^\infty(-\lambda ^2, -1; L^p) } (-t)^{ \frac{3-p}{2}}.
\end{align*}
Whence,  \eqref{2.4a}.

\hspace{0.5cm}
As a consequence of \eqref{2.4a} along with  Calder\'on-Zygmund's estimate get 
\begin{equation}
 \| \pi  (t)\|^{ \frac{p}{2}}_{ \frac{p}{2}} \le C \| u \|^p_{ L^\infty(-\lambda ^2, -1; L^p)} (-t)^{ \frac{3-p}{2}}  \quad  \forall\, t\in (-\infty, 0)
\label{2.4b}
\end{equation}
for a positive constant $ C$ depending only on $ p$ and $ \lambda $. In particular, from \eqref{2.4a} and \eqref{2.4b} respectively it follows that 
for all $ 2 \le q < \frac{2p}{p-3}$ and for all $ 0<R<+\infty$
\begin{equation}
u\in L^q(-R^2, 0; L^p(\R^{3})),\quad  \pi \in L^{ \frac{q}{2}}(-R^2,0; L^{ \frac{p}{2}}(\R^{3})),
\label{2.4bb}
\end{equation}  
 together with  the estimate
 \begin{align}
 \| u\|_{ L^q(-R^2, 0; L^p)}+  \| \pi \|^{ 1/2}_{ L^{ q/2}(-R^2, 0; L^{ p/2})} \le CR ^{ \frac{2}{q} - \frac{p-3}{p}} \| u \|_{ L^\infty(-\lambda ^2, -1; L^p)},
 \label{2.4bb1}
  \end{align} 
  where the constant $ C>0$ depends only on $ p$ and $ \lambda $. 
  
  \vspace{0.3cm}
{\it \underline{2. Local energy inequality}}: Let $ 0<r<R<+\infty$. Since $ u\in  C^\infty(Q)$ we get for all $ \phi \in C^{\infty}_{\rm c}
(B(0,R)\times (-R^2, 0])$ and for all $ t\in (-R^2, -r^2)$  the following   local energy equality 
\begin{align}
&\frac{1}{2} \intl_{ \R^{3}} | u(t)|^2  \phi^2 dx + \intl_{-R^2}^{t} \intl_{\R^{3}} | \nabla u|^2 \phi^2 dxds
\cr
&\qquad = \frac{1}{2}\intl_{-R^2}^{t} \intl_{\R^{3}} 
| u|^2 (\partial _t + \Delta ) \phi^2   dxds
+ \frac{1}{2} 
\intl_{-R^2}^{t} \intl_{\R^{3}} | u|^2 u\cdot \nabla \phi^2 dxds
\cr
&\qquad \qquad + \intl_{-R^2}^{t} \intl_{\R^{3}} \pi  u\cdot \nabla \phi ^2dxds
\cr
&\qquad = I +II+III.
\label{2.4f}
\end{align}

\hspace{0.5cm}
In \eqref{2.4f} we now take  $ \phi \in C^{\infty}(\R^{4})$ such that $ 0 \le \phi \le 1$ in $ \R^{4}$, 
$ \phi \equiv 1 $ on $ Q(0, R)$, $ \phi \equiv 0$ in $ Q  \setminus Q(0, 2R)$ and 
\[
| \nabla \phi | \le CR^{ -1}, \quad  | \partial _t \phi | + | \nabla ^2 \phi | \le CR^{ -2}. 
\]

\vspace{0.5cm}
\hspace{0.5cm}  
Firstly, employing H\"older's inequality togehter with  \eqref{2.4bb1}, we find
\[
I \le CR^{ -2} \| u\|^2_{2,  Q(0, 2R)} \le C R^{ 3- \frac{6}{p}- \frac{4}{q}} \| u\|_{ L^q(-4R^2, 0; L ^p)}^2 
\le C R  \| u \|^2_{ L^\infty(-\lambda ^2, -1; L^p)}. 
\] 
In particular, we have obtained the inequality 
\begin{equation}
 \| u\|^2_{2,  Q(0, 2R)} \le  C R^3  \| u \|^2_{ L^\infty(-\lambda ^2, -1; L^p)}. 
\label{2.4f1}
\end{equation}

Secondly, by the aid of H\"older's inequality  along with \eqref{2.4bb1} we estimate 
\begin{align*}
II &\le CR^{ -1} \intl_{Q(0, 2R)} | u|\phi  | u|^2 dxds   \le CR^{ -1}\| u \phi \|_{ L^{ \frac{q}{q-2}}(-4R^2, -r^2; 
L^{ \frac{p}{p-2}})}
\| u\|_{ L^q(-4R^2, 0; L ^p)}^2
\\
&\le C R ^{ \frac{4}{q} - \frac{2(p-3)}{p}-1} \| u \phi \|_{ L^{ \frac{q}{q-2}}(-4R^2, -r^2; 
L^{ \frac{p}{p-2}})}  \| u \|^2_{ L^\infty(-\lambda ^2, -1; L^p)}.
\end{align*}
In case $ p \ge 4$, having $ \frac{p}{p-2} \le 2$, with the help of Jensen's inequality we estimate 
\[
R ^{ \frac{4}{q} - \frac{2(p-3)}{p}-1}\| u \phi \|_{ L^{ \frac{q}{q-2}}(-4R^2, -r^2; 
L^{ \frac{p}{p-2}})} \le C R^{ \frac{1}{2}}\| u \phi \|_{ L^{\infty}(-4R^2, -r^2; 
L^{ 2})}
\] 
and by Young's inequality it follows that 
\[
II \le  \frac{1}{16} \| u \phi \|^2_{ L^{\infty}(-4R^2, -r^2; 
L^{ 2})} + C R \| u \|^4_{ L^\infty(-\lambda ^2, -1; L^p)}.
\]
In case $ 3 \le p < 4$ we choose  $ q = \frac{8p}{7p-12} $. As it readily seen  that   $ 2 < q < \frac{2p}{p-3}$ and $ \frac{2}{ \frac{q}{q-2}} + \frac{3}{ \frac{p}{p-2}}= \frac{3}{2}$.  Thus, by Sobolev's embedding theorem, H\"older's inequality, and  \eqref{2.4f1}, we obtain 
\begin{align*}
&\| u \phi \|_{ L^{ \frac{q}{q-2}}(-4R^2, -r^2; 
L^{ \frac{p}{p-2}})} 
\\
&\qquad \le C \| u \phi \|_{ L^{\infty}(-4R^2, -r^2; 
L^{ 2})} + \| \nabla u \phi\|_{ L^2(-4R^2, - r^2; L^2)} + C R^{ -1}\| u \|_{2, Q(0,2R)}  
\\
&\qquad \le C \| u \phi \|_{ L^{\infty}(-4R^2, -r^2; 
L^{ 2})} + \| \nabla u \phi\|_{ L^2(-4R^2, - r^2; L^2)} + C R^{ \frac{1}{2}}\| u \|_{ L^\infty(-\lambda ^2, -1; L^p)}.  
\end{align*}
Observing that $ \frac{4}{q} - \frac{2(p-3)}{p}-1= \frac{7p-12}{2p} - \frac{4p-12}{2p}-1 = \frac{1}{2}$, and  applying   Young's inequality, we find 
\begin{align*}
II &\le C R ^{ \frac{1}{2}} 
\Big(\| u \phi \|_{ L^{\infty}(-4R^2, -r^2; 
L^{ 2})} + \| \nabla u \phi\|_{ L^2(-4R^2, - r^2; L^2)}\Big)  \| u \|^2_{ L^\infty(-\lambda ^2, -1; L^p)} 
\\
&\qquad \qquad + C R\| u \|^3_{ L^\infty(-\lambda ^2, -1; L^p)}.  
\\
&\le 
\frac{1}{16} \| u \phi \|^2_{ L^{\infty}(-4R^2, -r^2; 
L^{ 2})} + \frac{1}{8}\| \nabla u \phi \|^2_{ L^2(-4R^2, - r^2; L^2)} 
\\
&\qquad \qquad \qquad + CR\Big(\| u\|^2_{ L^q(-4R^2, 0; L^p)} + 
\| u\|^4_{ L^q(-4R^2, 0; L^p)}\Big).  
\end{align*}
By an analogous reasoning using \eqref{2.4bb1},we get 
\begin{align*}
III &
\le CR^{ -1}\| u \phi \|_{ L^{ \frac{q}{q-2}}(-4R^2, -r^2; 
L^{ \frac{p}{p-2}})}
\| \pi \|_{ L^{ q/2}(-4R^2, 0; L ^{ p/2})}
\\
&\le C R ^{ \frac{4}{q} - \frac{2(p-3)}{p}-1} \| u \phi \|_{ L^{ \frac{q}{q-2}}(-4R^2, -r^2; 
L^{ \frac{p}{p-2}})}  \| u \|^2_{ L^\infty(-\lambda ^2, -1; L^p)}
\\
&\le \frac{1}{16} \| u \phi \|^2_{ L^{\infty}(-4R^2, -r^2; 
L^{ 2})} + \frac{1}{8}\| \nabla u \phi \|^2_{ L^2(-4R^2, - r^2; L^2)} 
\\
&\qquad \qquad \qquad + CR\Big(\| u\|^2_{ L^q(-4R^2, 0; L^p)} + 
\| u\|^4_{ L^q(-4R^2, 0; L^p)}\Big).  
\end{align*}
Inserting the above estimates of $ I, II$ and $ III$ into \eqref{2.4f}, and taking the supremum over $ (-R^2, -r^2)$ with respect to time, we arrive at 
\begin{align}
&\| u \phi \|^2_{ L^{\infty}(-4R^2, -r^2; 
L^{ 2})} + \| \nabla u \phi \|^2_{ L^2(-4R^2, - r^2; L^2)} 
\cr
&\qquad \qquad \le CR\Big(\| u \|^2_{ L^\infty(-\lambda ^2, -1; L^p)} + \| u \|^4_{ L^\infty(-\lambda ^2, -1; L^p)}\Big).
\label{2.4g}
\end{align}
Thus, by means of the lower semi-continuity of the norm,  letting $ r \rightarrow 0$ in \eqref{2.4g},  it follows that 
\begin{align}
&\| u\|^2_{ L^{\infty}(-R^2, 0; 
L^{ 2}(B(0,R)))} + \| \nabla u  \|^2_{2, Q(0,R)} 
\le  
CR\Big(\| u \|^2_{ L^\infty(-\lambda ^2, -1; L^p)} + \| u \|^4_{ L^\infty(-\lambda ^2, -1; L^p)}\Big),
\label{2.4c}
\end{align}
with a constant $ C>0$ depending only on $ p$ and $ \lambda $.  

\vspace{0.3cm}
{\it \underline{3. Serrin type estimate in terms of  weighted norm}}: Our next aim is to prove that  
\begin{equation}
 \intl_{-\infty}^{0}  \intl_{B(0,R)  \setminus B(0, R^{ -1})} | u |^p dx (-t)^{ \frac{p-5}{2}} dt <
 C \log R \| u\|^p_{ L^\infty(-\lambda ^2, -1; L^p)}\quad \forall\,1<R<+\infty,
\label{2.7a}
\end{equation}   
where the constant $ C>0$ depends only on $ p$ and $ \lambda $. 

\hspace{0.5cm}
In the proof below we use the following notation  
\[
D_k := B(0, \lambda ^{ -k+1}) \setminus B(0, \lambda ^{ -k}),  \quad k\in \Z. 
\]
Clearly, $ u\in L^p(\R^{3}\times (-\lambda ^2, -1))$ implies 
\begin{equation}
 \sum_{k=-\infty}^{\infty}  \intl_{-\lambda ^2}^{-1} \intl_{D_k} | u |^p dxdt  = \intl_{-\lambda ^2}^{-1}\intl_{ \R^{3}} | u |^p dxdt  
\le (\lambda ^2-1) \| u\|^p_{ L^\infty(-\lambda ^2, -1; L^p)}.
\label{2.3}
\end{equation}
By using the transformation formula of the Lebesgue integral we find 
\begin{align}
\intl_{-\lambda ^2}^{-1}\intl_{D_k} | u |^p dxdt &= \lambda ^{ (p-5)k}  \intl_{- \lambda ^{2k+2 }}^{-\lambda ^{ 2k}}   
\intl_{B(0, \lambda )  \setminus B(0,1)} | u_{\lambda^{ -k}} |^p dxdt 
\cr
& = 
\lambda ^{ (p-5)k} \intl_{- \lambda ^{2k+2 }}^{-\lambda ^{ 2k}}   \intl_{B(0, \lambda )  \setminus B(0,1)} | u |^p dxdt
\cr
&\ge \min\{1, \lambda^{5-p} \} \intl_{- \lambda ^{2k+2 }}^{-\lambda ^{ 2k}}   \intl_{B(0, \lambda )  \setminus B(0,1)}| u |^p (-t)^{ \frac{p-5}{2}}dxdt.
\label{2.4}
\end{align}

We now perform the sum over $ k\in \Z$ on both sides of \eqref{2.4}, which together with \eqref{2.3} 
gives 
\begin{align}
 &\intl_{-\infty}^{0}  \intl_{B(0, \lambda )  \setminus B(0,1) } | u |^p dx (-t)^{ \frac{p-5}{2}}  dt 
 \cr
 &\qquad \le  
(\min\{1, \lambda^{5-p} \})^{ -1} \sum_{k=-\infty}^{\infty}\intl_{-\lambda ^2}^{-1} \intl_{D _k} | u |^p dx  dt
\le C (\lambda -1) \| u\|^p_{ L^\infty(-\lambda ^2, -1; L^p)}.
\label{2.5}
\end{align}
Due to discrete self-similarity of $ u $ we get from \eqref{2.5} for all $ k\in \N$
\begin{equation}
 \intl_{-\infty}^{0}  \intl_{D_k} | u |^p dx (-t)^{ \frac{p-5}{2} } dt 
 = \intl_{-\infty}^{0}  \intl_{B(0, \lambda )  \setminus B(0,1) } | u |^p dx (-t)^{ \frac{p-5}{2}}  dt<
 C (\lambda -1) \| u\|^p_{ L^\infty(-\lambda ^2, -1; L^p)}.
 \label{2.6}
\end{equation}   
Given $ 1 < R< +\infty$ we may choose $ N \in \N$ such that $ \lambda ^N < R \le  \lambda ^{ N+1}$. 
Summation of \eqref{2.6} over $ k = -N+1, \ldots, N+1$ yields  
\begin{equation}
 \intl_{-\infty}^{0}  \intl_{B_{ \lambda ^{ N+1}}  \setminus B_{ \lambda ^{ -N}}} | u |^p dx (-t)^{ \frac{p-5}{2}} dt <
 CN(\lambda -1)  \| u\|^p_{ L^\infty(-\lambda ^2, -1; L^p)}\quad \forall\,N\in \N.
\label{2.7}
\end{equation}  
Since $ N (\lambda -1)\le C (\lambda -1)\frac{\log R}{\log \lambda }$, we get  \eqref{2.7a}.  We wish to remark that the constant 
in \eqref{2.7a} stays bounded as $\lambda \rightarrow 1 $.  

\vspace{0.2cm}  
{\it \underline{4. Regularity in $ \overline{Q}   \setminus \{ 0,0\}$}}: It suffices  to show that every point 
$z_0= (x_0, 0) \not= (0,0)$ 
is a regular point. Let us assume that    $z_0= (x_0, 0) \not= (0,0)$  is not a regular point, i.e. $ u$ not bounded in any neighborhood of $ z_0$. Appealing to \cite[Theorem\,5.1]{Wolf2015c},   
there exists {\color{blue} an} absolute  constant  $ \var >0$ such that 
\begin{equation}
r^{ -2}  \intl_{Q(z_0, r)} | u  |^3  dx dt \ge \var^3\quad  \forall 0< r < +\infty. 
\label{2.8}
\end{equation}
Otherwise, $ u$ is bounded in a neighborhood of $ z_0$. 
We also wish to emphasize that the above $ \var $ condition can be seen as an improvement of Scheffer's criterion 
(cf.  \cite{Sch1977}, and \cite[Theorem\,15.3]{Rob2016}, \cite[Lemma\,6.1]{Ser2015}), which includes the $ L^{ 3/2}(Q(z_0, r))$ norm of the pressure.  Although the above criterion has been  proved in  \cite[Theorem\,5.1]{Wolf2015c} for local suitable weak solution it still remain true in our case. Indeed, it not difficult to  check that if $ (u,p)$ is a  suitable weak solution    to the Navier-Stokes equations, then 
$ u$ is   a local suitable weak solution in the sense of \cite[Definition\,3.1]{Wolf2015c}.   For readers convenience a detailed proof of this claim is presented in the appendix of this paper.

\hspace{0.5cm}
Define $ \rho := \frac{1}{2}\min\{ | x_0|, 1\}$.  Then $Q(z_0, \rho ) \subset   B(0, 3\rho )  \setminus B(0, \rho ) \times (-\infty, 0)$. 
In view of \eqref{2.7a} we have  
\begin{equation}
 \intl_{Q(z_0, \rho )} | u|^p   dx  (-t)^{ \frac{p-5}{2}}   dt <+\infty. 
\label{2.9b}
\end{equation}
Let $ \{\rho _k\}$ be a sequence in $ (0, \rho )$  such that  $ \rho _k \rightarrow 0$ as $ k \rightarrow +\infty$. 
Then \eqref{2.9b} implies that 
\begin{equation}
 \intl_{Q(z_0, \rho _k)} | u|^p   dx  (-t)^{ \frac{p-5}{2}}   dt \rightarrow 0\quad  \text{as}\quad 
 k \rightarrow +\infty. 
\label{2.9}
\end{equation}

\hspace{0.5cm}
We  now define 
\begin{align*}
v_k (y,s) &= \rho _k u (x_0+ \rho _k y, \rho _k^2 s),\quad   
\\
\pi _k (y,s) &= \rho^2 _k \pi  (x_0+ \rho _k y, \rho _k^2 s),\quad  (y, s)\in \R^{3}\times (-\infty, 0).  
\end{align*}
 Then thanks to scaling invariance,   \eqref{2.8} yields 
\begin{equation}
 r^{ -2} \intl_{Q(0,r)} | v_k  |^3  dy ds \ge \var^3 \quad  \forall 0< r < +\infty. 
\label{2.10}
\end{equation}
On the other hand,  rescaling   \eqref{2.9} leads to  
 \begin{equation}
 \intl_{Q(0,1)} | v_k|^p   dy  (-s)^{ \frac{p-5}{2}}   ds \rightarrow 0\quad  \text{as}\quad 
 k \rightarrow +\infty. 
\label{2.11}
\end{equation}  
By means of Riesz-Fischer's theorem, eventually passing to a subsequence, from \eqref{2.11} we deduce that 
\begin{equation}
 v_k \rightarrow 0  \quad  \text{{\it a.e. in}}\quad  Q(0,1)\quad  \text{{\it as}}\quad  k \rightarrow +\infty.  
\label{2.11a}
\end{equation}

Furthermore, observing \eqref{2.4a} and \eqref{2.4b}, we infer that for all $ s\in (-\infty, 0)$ 
\begin{align}
\| v_k (s)\|_p^p &= \rho^{ p-3}_k \| u (\rho _k^2 s)\|_{ p} \le  
\rho^{ p-3}_k\lambda ^{ p-3}\| u \|^p_{ L^\infty(-\lambda ^2, -1; L^p) }(-\rho _k^2 s)^{ \frac{3-p}{2}}
\cr
& \le  \lambda ^{ p-3} \| u \|^p_{ L^\infty(-\lambda ^2, -1; L^p) } s^{ \frac{3-p}{2}}.
\label{2.12}
\end{align}
 In particular, for every  $ 2< q < \frac{2p}{p-3}$ and $ 0<R<+\infty$ we get 
the estimate 
\begin{equation}
 \| v_k\|_{ L^q(-R^2, 0; L^p)}+  \| \pi_k \|^{ 1/2}_{ L^{ q/2}(-R^2, 0; L^{ p/2})} \le C R ^{ \frac{2}{q} - \frac{p-3}{p}}\| u \|_{ L^\infty(-\lambda ^2, -1; L^p)} 
\label{2.14}
\end{equation}
with a constant $ C>0$ depending only on $ p$ and $ \lambda $. 
Since $ (v _k, \pi _k)$ is a solution to the Navier-Stokes equations, using the same argument as we have used in the proof of 
\eqref{2.4c},  from \eqref{2.14} we get for all $ 0<R<+\infty$ the estimate for the local energy 
\begin{align}
&\| v_k\|^2_{ L^{\infty}(-R^2, 0; 
L^{ 2}(B(0,R)))} + \| \nabla v_k  \|^2_{2, Q(0,R)} 
\cr
&\qquad \le  
CR\Big(\| u \|^2_{ L^\infty(-\lambda ^2, -1; L^p)} + \| u \|^4_{ L^\infty(-\lambda ^2, -1; L^p)}\Big).
\label{2.15}
\end{align}
By virtue of Sobolev's embedding theorem we see that $ \{ v_k\}$ is bounded in $ L^2(-R^2,0; L^6(B(0,R)))$, and thus by an interpolation 
argument it follows that $ \{ v_k\}$ is bounded in $ L^{ 10/3}(Q(0,R))$. Observing \eqref{2.11a},  we are now in a position to apply 
Vitali's comvergence theorem to conclude that 
\begin{equation}
 v _k \rightarrow 0   \quad  \text{{\it strongly in}}\quad  L^3(Q(0,1))\quad  \text{{\it as}}\quad  k \rightarrow +\infty. 
\label{2.18c}
\end{equation}

However, this contradicts \eqref{2.10}. 
Thus, we conclude that $u$ is regular on $\overline{Q}\setminus \{ (0,0)\}$.  We also wish to remark that from the convergence 
property \eqref{2.18c}  and the uniqueness of the limit after returning to the function $ u$, we obtain  
\[
\lim_{r \to 0} r^{ -2} \intl_{Q(z_0,r)} | u|^3 dxdt =0. 
\]
Then thanks to  \cite[Theorem\, 1.1]{gus} we  infer that  $ z_0$ is a regular point.

\vspace{0.2cm}  
 \underline{{\it 5. Proof of  (\ref{1.6a})}} .
 According to the step 4., where we have shown that $ u$ is bounded in any set $ \R^{3} \times (-\infty,0)  \setminus Q(0,r)$ 
 it holds  
\begin{equation}
|  u (x,t)| \le C \quad \qquad \forall\, (x,t) \in \overline{Q(0,\lambda )} \setminus Q(0,1). 
\label{2.19}
\end{equation}
Now let $ (x, t)\in \overline{Q}   \setminus \{(0,0)\}$. Then there exists $ k\in \Z$ such that 
\[
(x, t) \in \overline{Q(0, \lambda ^{ k+1})}  \setminus Q(0, \lambda ^k).
\]
Thus, $ (\lambda^{ -k} x, \lambda ^{ -2k} t) \in \overline{Q(0,\lambda )} \setminus Q(0,1).$ 
In view of \eqref{1.2} and \eqref{2.19}  it is readily seen that  
\[
| u(x,t) | = | u_{ \lambda ^{ -k}}(x,t)| = \lambda^{-k} | u(\lambda^{ -k} x, \lambda ^{ -2k} t)|\le  \lambda^{-k} C.   
\] 
As $\,\sqrt{-t} + | x| \le  2\max\{ | x|, \,\sqrt{-t}\}\le 2 \lambda^{ k+1}$, from the inequality above it follows that 
\[
| u(x,t) | \le \frac{2C\lambda }{\,\sqrt{-t}+| x|}.   
\]
This completes the proof of the theorem.  \hfill \Beweisende

\section{Proof of Theorem\,\ref{thm1.4}}
\label{sec:-3}
\setcounter{secnum}{\value{section} \setcounter{equation}{0}
\renewcommand{\theequation}{\mbox{\arabic{secnum}.\arabic{equation}}}}

Let $ u\in C((-\infty, 0); L^p(\R^{3}))\cap C^\infty (\Bbb R^3)$ be a smooth  solution to the Navier-Stokes equations satisfying  with a constant $ C_{  \ast}>0$ 
the inequality
 \begin{equation}
| u(x,t)| \le 
\frac{C_{ \ast}}{  \,\sqrt{-t}+ | x|}\quad  \forall\,(x,t)\in Q.
 \label{3.4}
 \end{equation}
Thus, by using the regularity theory of the Navier-Stokes equations we infer for all $ l\in \N$
\begin{equation}
| \nabla^l u(x,t) | \le   \frac{C_l}{  (-t)^{ \frac{1+l}{2}} + | x|^{ 1+l}} .
\label{3.6}
\end{equation} 
We denote by  $ \omega = \nabla \times u$ the vorticity of $ u$. From the vorticity equation and \eqref{3.6} 
 we deduce that 
\begin{equation}
| \partial _t\omega | \le \frac{C}{  (-t)^{2} + | x|^{4}},\quad  | \nabla \omega| 
\le \frac{C}{(-t)^{ \frac{3}{2}}+ | x|^3}.
\label{3.7}
\end{equation}

{\it 1. Condition for non trivial DSS functions}.  Thanks to \eqref{3.4} we get for $-R^2 \le  t< 0$,  $0<R<+\infty $
\[
\| u(t)\|^2_{ 2, B(0,R)} \le \intl_{B(0,R)} \frac{C}{-t + | x|^2} dx \le CR. 
\]
 Furthermore appealing to \eqref{3.6} with $ l=1$, we obtain for all $ k\in \N$
  \[
\| \nabla u\|^2_{ 2, Q(0, 2^{ -k+1}R )  \setminus Q(0, 2^{ -k}R) } \le \intl_{Q(0, 2^{ -k+1}R )  \setminus Q(0, 2^{ -k}R)} \frac{C}{(-t)^2 + | x|^4} dx \le C 2^{ -k} R. 
\]
Summation over $ k\in \N$ yields
\[
\| \nabla u\|^2_{ 2, Q(0,R) } \le  CR. 
\]
From the two estimates above we deduce that 
\begin{equation}
\| u\|^2_{ V^2(Q(0,R))} = \| u\|^2_{ L^\infty(-R^2, 0; L^2(B(0,R)))} + \| \nabla u\|^2_{2, Q(0,R)} \le CR. 
\label{3.4b}
\end{equation}

According to \cite[Theorem\,5.1]{Wolf2015c} there exists an absolute number $ \var >0$ 
such that if 
\begin{equation}
\intl_{Q(0,1)} | u|^3 dxdt \le  \var ^3,
\label{3.7a}
\end{equation}
then $ u$ is bounded in $ Q(0,1/2)$. 

\hspace{0.5cm}
We now assume that $ u$ is $ \lambda $-DSS $ (\lambda >1)$ and  non trivial.   Then we must have 
\begin{equation}
\intl_{Q(0,1)} | u|^3 dxdt >\var ^3.  
\label{3.7d}
\end{equation}
Otherwise, we get for every $ (x,t)\in Q$ 
\[
| u(x, t) | = \lambda^{ -k} | u(\lambda ^{ -k}, \lambda ^{ -2k}t)| \le C\lambda ^{ -k} \rightarrow 0
\]
as $ k \rightarrow +\infty$. 

\vspace{0.5cm}  
{\it 2. Indirect argument}.  We now assume the assertion of the theorem is not true. Then there exists a sequence 
$ \lambda _j \in (1, +\infty)$ with $ \lambda_j \rightarrow 1 $ as $ j \rightarrow +\infty$, 
and a sequence of non trivial $ \lambda _j$-DSS solutions $ u^j$ to the Navier-Stokes equations in $ Q$ satisfying  
the condtion \eqref{1.7}   for some constant $ C_{ \ast}>0$. Hence, by step 1 it follows that 
\begin{equation}
\intl_{Q(0,1)} | u^j|^3 dxdt >\var ^3\quad  \forall\,j\in \N.  
\label{3.7dd}
\end{equation}
Observing \eqref{3.4b}, by an reflexivity argument together with Banach-Alaoglu's theorem 
(eventually passing to a subsequence), we get a function $ u \in V^2_{ loc, \sigma } ( \overline{Q} )$  such that for all $ 0<R<+\infty$
\begin{align*}
\nabla  u^j &\rightarrow \nabla u  \quad  \text{{\it weakly in}}\quad  L^2(Q(0,R))\quad  \text{{\it as}}\quad  j \rightarrow +\infty 
\\
u^j &\rightarrow  u  \quad  \text{{\it weakly-$ \ast$ in}}\quad  L^\infty(-R^2,0; 
L^2(B(0,R))\quad  \text{{\it as}}\quad  j \rightarrow +\infty. 
\end{align*}
In order to verify the compactness with respect to the $ L^3(Q(0,R))$ norm we need  a priori bound for the pressure $ \pi ^j$. This can be done by decomposing $ \pi ^j$ into the sum $ \pi ^j_1+ \pi ^j_2$, where $ \pi _1^j$ and $ \pi _2^j$  is given by 
\[
-\Delta  \pi _1^j = \nabla \cdot \nabla \cdot ( u \otimes u \chi _{ | x| \le 1}) ,  \quad  -\Delta  \pi _2^j = \nabla \cdot \nabla \cdot  
(u \otimes u \chi _{ | x| >1 }).  
\]
Then by Calderon-Zygmund inequality together with \eqref{3.4}  we obtain for $ t\in (-\infty, 0)$
\begin{align*}
\| \pi _1^j(t)\|^{ 5/4}_{ L^{ 5/4}} &\le c\intl_{B(0,1)} | u^j|^{ 5/2} dx \le c C_{ \ast}^{ 5/2} \intl_{B(0,1)} |x|^{- 5/2} dx \le c C_{ \ast}^{ 5/2}, 
\\
\| \pi _2^j(t)\|^2_{ L^{ 2}} &\le c\intl_{B(0,1)^c} | u^j|^{ 4} dx \le c C^{ 4}_{ \ast} \intl_{B(0,1)^c} |x|^{- 4} dx \le c C_{ \ast}^{ 4}. 
\end{align*}
This shows that for all $ 0<R<+\infty$ we have the bound 
\begin{equation}
\| \pi ^j\|_{L^{ 5/4}(0, R)} \le c C_{ \ast}^2. 
\label{3.3e}
\end{equation}

 By means of compactness due to Aubin-Lions lemma we obtain  for all $ 0<R<+\infty$
\[
 u^j \rightarrow u  \quad  \text{{\it in}}\quad  L^3(Q(0,R))\quad  \text{{\it as}}\quad  j \rightarrow +\infty.
\]
With the help of the above convergence properties we infer that $ u$ is a local weak solution to the Navier-Stokes equations. In particular, from \eqref{3.7dd} we deduce that
 \begin{equation}
\intl_{Q(0,1)} | u|^3 dxdt \ge \var ^3\quad  \forall\,j\in \N.  
\label{3.7ddd}
\end{equation}

\hspace{0.5cm}
Let us now prove that $ u$ is backward self-similar. 
We set $ \omega ^j = \nabla \times u^j$, $ j\in \N$ and $ \omega =\nabla \times u$. 
Let $ \overline{Q(z_0, r)} \subset \overline{Q}   \setminus \{ (0,0)\}$, $ 0< r<+\infty$. 
According to  \eqref{3.7}  $ | \partial _t \omega ^j|$ and $ | \nabla \omega ^j |$ are uniformly 
bounded on $ \overline{Q(z_0, r)}$.  Using Arzel\`a-Ascoli's theorem, eventually passing to a subsequence, we get 
\begin{equation}
 \omega ^j  \rightarrow \omega   \quad  \text{{\it uniformly on $ \overline{Q(z_0, r)}$}}\quad \text{{\it as}}\quad  j  \rightarrow +\infty.  
\label{3.8}
\end{equation}

\hspace{0.5cm}
 Let $ 1< \mu < +\infty$ be arbitrarily chosen, but fixed. 
 Clearly, there exists a unique  $ k_j= k(\lambda_j , \mu )\in \N$ such that 
\begin{equation}
\lambda_j ^{ k_j-1} \le \mu < \lambda_j ^{ k_j}. 
\label{3.4a}
\end{equation}
We observe that $0 \le \lambda ^{ k_j}_j - \mu  \le  \lambda_j ^{ k_j} - \lambda ^{ k_j-1}_j = \lambda ^{ k_j-1}_j (\lambda_j -1) \le \mu (\lambda_j -1) \rightarrow 0$ as $ j  \rightarrow +\infty$, which shows that 
\begin{equation}
\lambda_j ^{ k_j} \rightarrow  \mu \quad  \text{ {\it as}}\quad  j \rightarrow +\infty. 
\label{3.5a}
\end{equation}

Let $z= (x,t)\in Q$. 
Note that for $ j\in \N$ taken sufficiently large, we have 
 $ | x - \lambda _j ^{ k_j} x | = (\lambda _j^{ k_j}-1) | x| \le \lambda _j 
\mu | x| \le \mu ^2 | x| $, and 
$ | t -\lambda _j^{ 2k_j} t|^{ \frac{1}{2}} \le (\lambda _j^{ 2k_j}-1)^{ \frac{1}{2}} | t|^{  \frac{1}{2} } \le \mu ^2 \,\sqrt{-t}$. Setting $ R= \mu ^2\max\{ | x|, \,\sqrt{-t}\}$ we see that for sufficiently large $ j\in \N$
\[
z^{j}=(\lambda ^{ k_j}_j x, \lambda ^{ k_j}_j t) \in \overline{Q(z, R)} \subset \R^{3}\times (-\infty, 0]  \setminus \{ (0,0)\}. 
\]
In addition, it can been easily checked that $ (\mu x, \mu ^2 t) \in \overline{Q(z, R)} $.
Using triangular inequality and the fact that 
\[
\lambda_j^{ 2k_j} \omega ^{ j} (\lambda_j^{ k_j} x, \lambda_j^{ 2k_j} t) = 
\omega^j (x, t), 
\]
 we get  
\begin{align*}
&| \mu^2 \omega (\mu x, \mu ^2 t) - \omega (x,  t)| 
\\
&\quad \le  | \mu^ 2\omega (\mu x, \mu ^2 t) - 
\mu^2 \omega ^j (\mu x, \mu ^2 t)|
+ | \mu^2 \omega ^j (\mu x, \mu ^2 t) -  \omega ^j (x,  t)| 
\\
& \qquad +| \omega ^j  (x,  t) - \omega (x,  t)|.
\end{align*}
Clearly, thanks to \eqref{3.8} the first term and the last term  on the right-hand side converges to zero as $ j  \rightarrow +\infty$. We only need to 
investigate the second term. In fact, by using the discrete self-similarity of each $ \omega ^j$ and triangular inequality, we find 
\begin{align*}
& | \mu^2 \omega ^j (\mu x, \mu ^2 t) -  \omega ^j (x,  t)| 
\\
&\quad  =| \mu^2 \omega ^j (\mu x, \mu ^2 t) -  \lambda_j^{ 2k_j} \omega ^j (\lambda^{ k_j}_j x,  \lambda_j ^{ 2k_j} t)|
\\
& \quad \le  ( \lambda_j^{ 2k_j}-\mu^2)|  \omega ^j (\mu x, \mu ^2 t)| 
+ \lambda_j^{ 2k_j} | \omega ^j (\mu x, \mu ^2 t) - \omega ^j (\lambda^{ k_j}_j x,  \lambda_j ^{ 2k_j} t)|
\\
& \quad \le  
  ( \lambda_j^{ 2k_j}-\mu^2)|  \omega ^j (\mu x, \mu ^2 t)| 
+  \lambda_j^{ 2k_j} | \omega ^j (\mu x, \mu ^2 t) - \omega (\mu  x, \mu ^2 t)|
\\
& \qquad +\lambda_j^{ 2k_j} | \omega  (\mu x, \mu ^2 t) - \omega  (\lambda^{ k_j}_j x,  \lambda_j ^{ 2k_j} t)| +
\lambda_j^{ 2k_j} | \omega (\lambda^{ k_j}_j x, \lambda^{ 2k_j}_j t) - \omega ^j (\lambda^{ k_j}_j x,  \lambda_j ^{ 2k_j} t)|.
\end{align*} 
It is readily seen that due to \eqref{3.8}  and \eqref{3.5a}  the first term tends to zero as $ j\rightarrow +\infty$, while by virtue of  \eqref{3.8}, \eqref{3.5a}
and the continuity of $ \omega $,  the second, third   and fourth   
 term tends to zero as $ j \rightarrow +\infty$.  Consequently, 
\[
\mu^2 \omega (\mu x, \mu ^2 t) = \omega (x,t). 
\] 
In particular, $ \nabla \times (u_\mu - u) =0 $. Due to $ \nabla \cdot (u_\mu - u)=0$ the function $ u_{ \mu }- u$ 
is  harmonic. Observing   \eqref{3.4},  the Liouville theorem for harmonic functions implies $ u_{ \mu }- u=0$. Hence, $ u$ is a backward self-similar solution to the Navier-Stokes equations 
fulfilling 
\begin{equation}
| u(x,t)| \le \frac{C}{\,\sqrt{-t}+ |x |},\quad  (x,t)\in \overline{Q}   \setminus \{ (0,0)\}.
\label{3.10}
\end{equation}
In particular, $ u$ satisfies  the local energy estimate \eqref{3.4b}. Thus,  we are in a position  to apply Tsai's result 
\cite[Theorem\,2]{Tsai1998}, to see that $ u$ is identical zero. However this  contradicts to \eqref{3.7ddd}. 
Since our assumption is false the assertion of the theorem must be true. 
 \hfill \Beweisende 

\section{Proof of Theorem\,\ref{thm1.5}}
\label{sec:-4}
\setcounter{secnum}{\value{section} \setcounter{equation}{0}
\renewcommand{\theequation}{\mbox{\arabic{secnum}.\arabic{equation}}}}
Although the proof is similar to the proof of Theorem 1.2 in \cite{cha1}, we write it in detail for reader's convenience.\\
\ \\
We rewrite $v$ in terms of the self-similar variables as
$$v(x,t)=\frac{1}{\sqrt{t_*-t}}V\left( \frac{x-x_*}{\sqrt{t_*-t}}, -\log (t_*-t) \right)
$$
for some $V\in C(-\infty, +\infty; L^p (\Bbb R^3)\cap C^\infty (\Bbb R^3))$,  and transform $(u, \pi)  \to (U, P)$  by the formula 
$$u(x,t)=\frac{1}{\sqrt{t_*-t}} U(y,s), \quad \pi (x,t)= \frac{1}{t_*-t} P (y,s), $$
where
$$ y= \frac{x-x_*}{\sqrt{t_*-t}},\quad s= -\log (t_*-t) 
$$
Then, we notice that the condition of discrete  self-similarity of the function $u(x,t)$  of  (\ref{1.2}) is equivalent to the time-periodicity of $U$, $V(\cdot, s)=V(\cdot, s+2\log \l )$ for all $s\in \Bbb R$, and  $(U, P)$ solves
\bb\label{ssn}
U_s +\frac12 U +\frac12 (y\cdot \nabla )U +(U\cdot \nabla ) U -\Delta U=-\nabla P, \quad \nabla \cdot  U=0,
\ee
and
the condition (\ref{as}) is transformed into
\bb\label{ss1}
\lim_{s\to \infty} \|U (\cdot, s)-V(\cdot, s)\|_{L^p (B_R (0))} =0 \quad \forall  R>0.
\ee
We also note that the discrete  self-similarity $\l u(\l x, \l^2 t)=u(x,t)$ is equivalent to the time periodicity
$$
U(\cdot , s)= U(\cdot , s+S_0 ), \quad S_0 :=2\log \l .
$$
Given $\xi \in C_0 ^\infty (0, S_0)$,  $\phi =(\phi_1, \phi_2, \phi_3) \in C_0 ^\infty (\Bbb R^3)$  with $\nabla \cdot \phi=0$ and $n\in \Bbb N$, we take $L^2 (\Bbb R^3 \times [n, n+S_0])$ inner
product  the first equation of (\ref{ssn}) by $\xi (\cdot-S_0 n) \phi $. Then, after integration by part we obtain
\bq\label{ss2}
\lefteqn{-\int_0 ^{S_0} \int_{\Bbb R^3} \xi_s (s)\phi (y) \cdot V(y, s+S_0 n) dyds -\int_0 ^{S_0} \int_{\Bbb R^3} \xi (s)\phi (y) \cdot V(y, s+S_0 n) dyds}\n \\
&&\quad -\frac12 \int_0 ^{S_0} \int_{\Bbb R^3} \xi(s)\phi (y) \cdot (y\cdot \nabla)V(y, s+S_0 n) dyds \n \\
&&\quad-\int_0 ^{S_0} \int_{\Bbb R^3} \xi(s)[ V(y, s+S_0 n)\cdot V(y, s+S_0 n)\cdot \nabla )\phi (y) ] dyds  \n \\
&&=\int_0 ^{S_0} \int_{\Bbb R^3} \xi (s)V(y, s+S_0 n) \cdot \Delta \phi (y)  dyds
\eq
Similarly from the second equation of (\ref{ssn}) we have
\bb \label{ss3}
\int_0 ^{S_0} \int_{\Bbb R^3} \xi (s)  V(y, s+S_0 n)\cdot \nabla\psi (y)dyds=0.
\ee
for all $\psi \in C_0 ^\infty (\Bbb R^3)$
Passing $n \to \infty$ in (\ref{ss2}) and (\ref{ss3}) and  recalling (\ref{ss1}),  we find that $V$ satisfies
$$
\int_0 ^{S_0} \int_{\Bbb R^3} \left\{ V_s +\frac12 V +\frac12 (y\cdot \nabla )V +(V\cdot \nabla )V-\Delta V\right\} \cdot \phi (y) \xi(s) dyds =0
$$
for all $\phi \in [C_0 ^\infty (\Bbb R^3)]^3$ with $\nabla \cdot \phi=0$ and $\xi \in C_0 ^\infty (0, S_0)$, and 
$$
\int_0 ^{S_0} \int_{\Bbb R^3} \xi (s) [\nabla \cdot  V] \psi (y)dyds=0.
$$
for all $\psi \in C_0 ^\infty (\Bbb R^3), \xi\in C_0 ^\infty (0, S_0). $ Therefore there exists $\bar{P}$ such that
\bb\label{ssn4}
V_s +\frac12 V +\frac12 (y\cdot \nabla )V +(V\cdot \nabla ) V -\Delta V=-\nabla \bar{P}, \quad \nabla \cdot  V=0,
\ee
Since $(v,\pi)$ given by
$$v(x,t)=\frac{1}{\sqrt{t_*-t}}  V(y,s), \quad \pi (x,t)= \frac{1}{t_*-t} \bar{P} (y,s), $$
is a discretely self-similar solution with the scaling parameter $\l \in (1, \l_*)$, applying Theorem \ref{thm1.4}, we find $v=V=0$, and 
(\ref{ss1}) reduces to
\bb\label{ss5}
\lim_{s\to \infty} \|U (\cdot, s)\|_{L^p (B_R (0))} =0 \quad \forall  R>0,
\ee
which can be written, in terms of the physical variables, as
\bb \label{ss6}
\lim_{t\to t_*} \left\{ (t_*-t)^{\frac{p-3}{2p}} \sup_{t<\tau<t_*} \|u(\cdot , \tau )\|_{L^p (B_{R \sqrt{ t_*-t}} (x_*))} \right\} =0.
\ee
Setting $R=1$, and $\sqrt{t_*-t}=r$ in (\ref{ss6}), we have
\bb \label{ss7}
\lim_{r\to0} \left\{ r^{\frac{p-3}{p}} \sup_{-r^2<\tau<0} \|u(\cdot , \tau )\|_{L^p (B_{r} (x_*))} \right\}=0.
\ee
Applying the regularity criterion by Seregin-\v{S}ver\'ak (cf.  \cite[Lemma\,3.3]{SerSve2002}), we are led to the fact that $z_*=(x_*, t_*)$ is a regular point.   \hfill \Beweisende 
\hspace{0.5cm}
$$\mbox{\bf Acknowledgements}$$
Chae was partially supported by NRF grant 2016R1A2B3011647, while Wolf has been supported by 
the German Research Foundation (DFG) through the project WO1988/1-1; 612414.

 \appendix
\section{Remark on the  notion of local suitable weak solutions}
\label{sec:-A}
\setcounter{secnum}{\value{section} \setcounter{equation}{0}
\renewcommand{\theequation}{\mbox{A.\arabic{equation}}}}

In this appendix we would like to clarify that any suitable weak solution to the Navier-Stokes equations  is a local suitable weak solution
in the sense of  \cite[Definition\,3.1]{Wolf2015c}. 

\hspace{0.5cm}
To this end, let $ \Omega \subset \R^{3}$ be a  domain and  $ 0< T<+\infty$.  By $ Q$ we denote the 
space time cylinder $ \Omega \times (0, T)$. We denote  $ V^{ 1,2}_\sigma (Q)$ the space of all vector functions  
$ L^\infty(0,T; L^2(\Omega ))\cap L^2(0,T; W^{1,\, 2}(\Omega ))$ fulfilling $ \nabla\cdot u =0 $ a.e. in $ Q$.  
We recall the following   notion of localized suitable weak solution to the Navier-Stokes equations,  which  is more general than 
the usual notion given by Scheffer in \cite{Sch1977}.  

\begin{defin}
Given $  f \in L^{2}(Q; \R^{3})$, a pair $ (u, p) \in V^{ 1,2}_\sigma (Q)\times L^{ 3/2}(Q)$
is called a {\it suitable weak solution} to the Navier-Stokes equations 
\begin{align}
\nabla \cdot u &=0 \quad  \text{ in}\quad Q, 
\label{A.1}
\\
\partial _t u + u\cdot \nabla u - \Delta u &= - \nabla p + f\quad  \text{ in}\quad Q,
\label{A.2}
\end{align}
if  \eqref{A.2} is satisfied in the sense of distributions, i.e. for every $ \varphi \in C^\infty(Q)$ it holds 
\begin{align}
\intl_{Q} - u \cdot \partial _t \varphi + (- u \otimes u  + \nabla u) : \nabla \varphi dxdt
 = \intl_{Q} p \nabla \cdot \varphi + f \cdot \varphi dxdt, 
\label{A.3}
\end{align}

 and if the the following local energy inequality holds true for a.e. 
$ 0< t< T$ and for all non negative $ \phi \in C^{\infty}(Q)$ with $ \supp(\phi ) \subset \Omega \times (0, T]$,     
\begin{align}
 &\intl_{\Omega } | u(t)|^ 2 \phi(t) dx + 2 \intl_{0}^{t} \intl_{\Omega}  | \nabla u|^2 \phi   dx  ds
\cr
&\qquad \le   \intl_{0}^{t} \intl_{\Omega}  | u|^2 (\partial _t + \Delta ) \phi  + (| u|^2 + 2p) u\cdot \nabla \phi  + 
2f\cdot u \phi  dx  ds. 
\label{A.4}
\end{align}

\end{defin}

\begin{lem}
\label{lemA.1} 
Let $ (u, p) \in V^{ 1,2}_\sigma (Q)\times L^{ 3/2}(Q)$ be a suitable weak solution to the Navier-Stokes equations \eqref{A.1}, \eqref{A.2}. Then $ u$ is a local suitable weak solution in the sense of  \cite[Definition\,3.1]{Wolf2015c}.    
\end{lem}

{\bf Proof}:  Let $ B \subset \Omega $  be a fixed ball. 
Since $ B $ is  bounded, we easily see that all terms in \eqref{A.2} besides $ \partial _t u$ belong to 
$ L^{ 3/2}(0, T; W^{-1,\, 3/2}(B ))$. Accordingly, $ u$ admits a  distributional time derivative in  $ L^{ 3/2}(0, T; 
W^{-1,\, 3/2}(B ))$. Let $ E^{ \ast}_B = \nabla \mathscr{ P}_B$ denote the local pressure projection introduced in \cite{Wolf2015c},  which is a projection in $ W^{-1,\, r}(B), 1 <r< +\infty$,  onto the closed subspace of functionals of the form $ \nabla q$ 
(for more details on the properties of $ E^{ \ast}_B$ we refer to  \cite{wol2}). 

\hspace{0.5cm}
Now, we define
\[
\nabla p_{ h, B}(t)= - E^{ \ast}_B (u(t)),\quad  \nabla p_{ 0, B} = E^{ \ast}_B (-  (u\cdot \nabla) u + \Delta u + f), 
\]
and set $ v_B = u + \nabla p_{ h, B}$.  We also wish to note that owing to $ \nabla \cdot u =0$  
the pressure $ p_{ h, B}(t)$ is harmonic in $ B$ for a.e. $ t\in (0,T)$.    

\hspace{0.5cm}
Since $ E^{ \ast}_B$ is a bounded spatial operator, it commutates with the  distributional time derivative. This  yields 
\[
E^{ \ast}_B(\partial _t u) = \partial _t (E^{ \ast}_B (u)) = - \partial _t \nabla  p_{ h, B}\quad  in \quad  B\times (0,T)
\]  
in the sense of distributions. Since the left-hand side belongs to $ L^{ 3/2}(0,T; W^{-1,\, 3/2}(B))$, we infer that 
$ \nabla p_{ h, B}$  admits a  distributional time derivative in     $ L^{ 3/2}(0, T; W^{-1,\, 3/2}(\Omega ))$. Thus taking into account 
that $ 	p_{ h, B}$ is harmonic with respect to the spatial variable, using the mean value property of harmonic functions together 
with Caccioppoli   inequality, it follows that $ \partial _t \nabla p_{ h, B} \in L^{ 3/2}(0, T; C^\infty(B) )$.   On the other hand, 
as $\nabla  p\in L^{3/2 }(0, T; W^{-1,\, 3/2}(B))$ we get 
\begin{align}
 \nabla p = E^{ \ast}_B (\nabla p) &= - E^{ \ast}_B(\partial _t u) +  E^{ \ast}_B (-  (u\cdot\nabla) u + \Delta u + f) 
 \cr
 &= \partial _t 
 \nabla p_{ h, B} + \nabla p_{ 0, B}
 \label{A.5}
\end{align}
in $ B\times (0,T)$ in the sense of distributions. 

\hspace{0.5cm}
Let $ 0< t< T$ be chosen    so that 
\begin{equation}
\frac{1}{h}  \intl_{t}^{t+h} u(s)  ds \rightarrow  u(t)\quad \text{ {\it in}}\quad  L^2(B)\quad  \text{ {\it as}}\quad  h \rightarrow 0^+.   
\label{A.6}
\end{equation}
Clearly,  by means of the   boundedness of the operator $ E^{ \ast}_B$ in $ L^2(B)$, the condition \eqref{A.6} implies   
\begin{equation}
\frac{1}{h}  \intl_{t}^{t+h} \nabla p_{ h, B}(s)  ds \rightarrow  \nabla p_{ h, B}(t)\quad \text{ {\it in}}\quad  L^2(B)\quad  \text{ {\it as}}\quad  h \rightarrow 0^+.   
\label{A.7}
\end{equation}

Let $ \phi \in C^{\infty}(Q) $  with $ \phi  \ge 0$ and $ \supp (\phi ) \subset B\times (0, T]$ be arbitrarily chosen. For 
$ 0< h < T-t$ by $ \eta _h \in C^{ 0,1}(\R)$ we denote the piecewise linear function such that $ \eta _h \equiv 1$ in $ (-\infty, t]$ 
and $ \eta _h= 0 $ in $ [t+h, T)$. Surely, $ \eta' = -\frac{1}{h} \chi _{ (t, t+h)}$.   
By a routine density argument it is readily seen that $ \varphi =2\nabla p_h \phi \eta _h $ is an admissible test function for \eqref{A.3}. 
Inserting this function into \eqref{A.3}, using the identity \eqref{A.5}, and applying integration by parts, we obtain 
\begin{align}
&\frac{2}{h} \intl_{t}^{t+h}  \intl_{B} u\cdot \nabla p_{ h,B} \phi  dxds -  
2\intl_{0}^{t+h}  \intl_{B}  u \cdot \nabla p_{ h, B}\partial _t \phi \eta _h + u \cdot \partial _t\nabla p_{ h, B} \phi \eta _hdxds
\cr
&\qquad \qquad  +2\intl_{0}^{t+h}  \intl_{B} (-u \otimes u  + \nabla u) : \nabla (\nabla p_{ h, B} \phi ) \eta _h dxds
\cr
&\qquad  = 2\intl_{0}^{t+h}  \intl_{B} p \nabla p_h \cdot \nabla \phi \eta _h + f \cdot \nabla p_h \phi \eta _h dxds
\cr
&\qquad  = - \frac{1}{h} \intl_{t}^{t+h}  \intl_{B} |  \nabla p_{ h,B}|^2 \phi  dxds + 
2\intl_{0}^{t+h}  \intl_{B} | \nabla  p_h|^2\partial _t\phi \eta _h dxds
\cr
&\qquad \qquad + 
\intl_{0}^{t+h}  \intl_{B} p_0 \nabla p_h \cdot \nabla \phi \eta _h + f \cdot \nabla p_h \phi \eta _h dxds.
\label{A.8}
\end{align}
Thanks to \eqref{A.6} and \eqref{A.7} we are in a position to pass $ h \rightarrow 0$ in both sides of \eqref{A.8}. This together with 
an elementary manipulation of the resultant identity we obtain 
\begin{align}
&  \intl_{B} (2u(t)\cdot \nabla p_{ h,B}(t) + | \nabla p_{ h, B}(t)|^2)\phi(t) dxds 
\cr
&\quad =  
\intl_{0}^{t}  \intl_{B}  (2u \cdot \nabla p_{ h, B}+ | \nabla p_{ h, B}|^2)\partial _t \phi dxds
\cr
&\quad \qquad  +2\intl_{0}^{t}  \intl_{B} (u \otimes u  - \nabla u) : \nabla (\nabla p_{ h, B} \phi ) dxds
\cr
&\quad\qquad  - 2\intl_{0}^{t}  \intl_{B} u \partial _t  p_{ h, B} \cdot  \nabla \phi dxds + 
2\intl_{0}^{t}  \intl_{B}  p_{ 0, B} \nabla p_{ h, B}\cdot  \nabla \phi  +  f \cdot \nabla p_h \phi  dxds. 
\label{A.9}
\end{align} 
Combining \eqref{A.9} and \eqref{A.4}, and once more appealing to \eqref{A.5},   we are led to 
\begin{align}
 &\intl_{\Omega } | v_B (t)|^ 2 \phi(t) dx + 2 \intl_{0}^{t} \intl_{\Omega}  | \nabla v_B |^2 \phi   dx  ds
 \cr
&\qquad \le   \intl_{0}^{t} \intl_{\Omega}  | v_B |^2 (\partial _t +\Delta)  \phi  + (| u|^2 u+ 2p_{ 0, B} v_B )\cdot \nabla \phi  + 
f\cdot u \phi  dx  ds
\cr
&\qquad \qquad  +2\intl_{0}^{t}  \intl_{B} u \otimes u : \nabla (\nabla p_{ h, B} \phi)  dxds.
\label{A.10}
\end{align}
Hence $ u$    is  a local suitable weak solution to \eqref{A.1}, \eqref{A.2} in the sense of 
 \cite[Definition\,3.1]{Wolf2015c}.  \hfill \Beweisende

 \begin{rem}
 By a slight modification of the above proof  it is readily seen that the statement of  Lemma\,\ref{lemA.1} remains valid 
  even if we replace $ f \in L^2(Q)$ by a more general right-hand side  $ f \in L^2(0,T; W^{-1,\, 2}(\Omega ))  + L^1(0,T; L^2(\Omega ))$. 
 \end{rem}

\end{document}